# A spatial local method for solving 2D and 3D advection-diffusion equations

## Huseyin Tunc, Murat Sari*

Department of Mathematics, Faculty of Arts and Science, Yildiz Technical University, Istanbul 34220, Turkey.

**Abstract**

In this study, an implicit-explicit local differential transform method (IELDTM) based on Taylor series representations is produced for solving 2D and 3D advection-diffusion equations. The parabolic advection-diffusion equations are reduced to the nonhomogeneous elliptic system of partial differential equations with the utilization of the Chebyshev spectral collocation approach in temporal variable. The IELDTM is constructed over 2D and 3D meshes using continuity equations of the neighbour representations with either explicit or implicit forms in related directions. The IELDTM is proven to have excellent convergence properties by experimentally illustrating both $h-$refinement and $p-$refinement outcomes. A distinctive feature of the IELDTM over existing numerical techniques is the optimization of the local spatial degrees of freedom. It has been proven that IELDTM provides more accurate results with far less degrees of freedom than the finite difference, finite element and spectral methods.

## 1. Introduction

Many quantities are encountered in various fields of science, such as mass, heat, energy, velocity, and concentration described in the advection-diffusion (AD) equation as the dependent variable. The (2+1) and (3+1) dimensional AD equations occur in the modelling of various real-life problems such as the tumour angiogenesis model, the tumour invasion model, heat transfer in draining film, dispersion of tracers in porous media, the spread of pollutants in rivers and stream, the dispersion of dissolved material in estuaries and coastal sea and so on [1]. When the advection is dominant to the diffusion in the equation, the exact series-based solutions mostly fail and thus diverge. It is not easy to find out analytical solutions in nonlinear problems even the problem is not advection-dominated. In such scenarios, effective numerical methods need to be constructed to observe the behaviour of the advection-diffusion mechanism.

---

*Corresponding author: E-mail: tnchsyn@gmail.com (Huseyin Tunc)

Various versions of the finite element method (FEM) [2-5], discontinuous Galerkin method (DGM) [6-7], finite difference method (FDM) [8-11], finite volume method (FVM) [12-13] and meshless method [14] have been presented for spatial parts of the AD equations in 2D and 3D. Among these numerical methods, finite element-based techniques have been the most preferred method for solving AD equations due to the geometric flexibility and weak formulation of the method [2-5]. Despite all these considerable advantages, the spatial degrees of freedom (*dof*) of the FEM increase rapidly with the use of the $p$ −refinement procedure in two and three space dimensions. Thus, the FEM in higher dimensions cannot take advantage of $p$ −refinement, and the $h$ −refinement procedure seems to be the only reasoned option to increase the accuracy of the results. The discontinuous Galerkin methods (DGM) have several advantages: high geometric flexibility, mass conservation properties, and shock-capturing nature [6-7]. Since the common nodes of the elements contribute two or more times in the DGM for solving (2+1) and (3+1) dimensional AD equations, both $h$ −refinement and $p$ −refinement procedures lead to a considerable amount of the degrees of freedom. In addition to these methods, the finite difference methods were widely used for solving the AD equations either in explicit or implicit forms [8-11]. The FDM takes advantage of $p$ −refinement procedure by increasing the nodes used in the stencils. The required degrees of freedom are not affected by order of the FDM formulation. However, the FDM has crucial restrictions, such as the need for a structured geometric domain and the necessity of using low order approximations around the boundaries. One more disadvantage of the FDM is the high-order formulations are not directly related to each other. Therefore, to increase the order of the method, all formulations must be redefined. All of the mentioned methods above for solving AD equations are multi-point numerical methods, i.e. the derivatives that occurred in the equations are approximated with the use of local differential operators. This group of methods uses the differential equations in various forms, such as integral forms as in the FEM or collocation forms as in the FDM, to obtain an algebraic system of equations. On the other hand, the Taylor series-based methods use differential equations to calculate higher-order derivatives with priori information such as the positions or the slopes of the solution functions. These higher-order derivatives are used to represent the approximate local or global solutions. The global semi-analytic Taylor series-based methods are known to be highly restricted for some specific problems [21]. The Taylor series-based multi-center techniques have been widely used to solve higher dimensional elliptic partial differential equations that occurred in elasticity problems [15-20]. The primary motivation of using local Taylor series expansions is the minimization of required degrees of freedom with the use of recurrence relation obtained from the differential equation itself. The nature of the elliptic PDEs is suitable for the multi-center Taylor series-based methods since one has information at all boundaries to be used in interior nodes.

Even the literature is quite rich for solving elliptic PDEs by Taylor series-based approaches, the derivation of a higher-order numerical method with optimized degrees of freedom based on the Taylor series for solving parabolic PDEs will fill a serious gap. In this study, an implicit-explicit local differential transform method (IELDTM) based on Taylor series expansions is derived for solving the (2+1) and (3+1) dimensional AD equations. Both the (2+1) and (3+1) dimensional AD equations are converted to the system of two- and three-dimensional nonhomogeneous elliptic boundary value problems with the use of the Chebyshev spectral collocation method (ChSCM) as a continuous time integrator. It is proven that the ChSCM provides excellent interpolation properties with less temporal *dof* for time integration of the AD equations. Both $h-$convergence and $p-$convergence of the present IELDTM are presented for 2D and 3D problems. It is shown that the present algorithm provides much more accurate results with fewer *dof* than the finite difference methods, the finite element methods, and the spectral methods for solving linear/nonlinear 2D equations and linear 3D equations.

## 2. IELDTM-ChSCM for (2+1)-dimensional AD equation

In this section, we introduce the IELDTM-ChSCM hybridization for solving the following (2+1) dimensional AD equation

$$u_t + V_x u_x + V_y u_y = D_x u_{xx} + D_y u_{yy}, \quad (x,y) \in A \text{ and } t \in [0, t_f] \tag{1}$$

with Dirichlet boundary conditions

$$u(a, y, t) = f_1(y, t) \text{ and } u(b, y, t) = f_2(y, t) \tag{2}$$

$$u(x, c, t) = h_1(x, t) \text{ and } u(x, d, t) = h_2(x, t) \tag{3}$$

and initial condition

$$u(x, y, 0) = g(x, y) \tag{4}$$

where $A = [a, b] \times [c, d]$, $V_x$ and $V_y$ are advection constants, $D_x$ and $D_y$ are diffusion constants and $f_1$, $f_2$, $h_1$, $h_2$ and $g$ are known smooth functions. The subscripts $x, y$, and $t$ represent differentiations with respect to spaces $x, y$, and time $t$, respectively. In the following subsection, (2+1) dimensional AD equation (1) will be reduced to a system of elliptic boundary value problems with the use of the Chebyshev spectral collocation method (ChSCM) in time.

**2.1. Temporal variation for (2+1) dimensional AD equation**

Let us approximate the time part of the $u(x, y, t)$ with $Nth$ order Chebyshev polynomials as follows [22]

$$u(x, y, t) = \sum_{n=0}^{N} \alpha_n c_n(x, y) \bar{T}_n(t) \tag{5}$$

where $c_n(x, y)$ represent the spatial parts of the $u(x, y, t)$, $\alpha_n = 1$ for interior points and $\alpha_0 = \alpha_N = 1/2$ for any time boundary point. In equation (5), $\bar{T}_n(t)$ is defined in terms of the first kind of Chebyshev polynomials $T_n(t)$ as

$$\bar{T}_n(t) = T_n\left(\frac{2t-t_f}{t_f}\right) = \cos\left(n\cos^{-1}\left(\frac{2t-t_f}{t_f}\right)\right). \tag{6}$$

We define the following restricted collocation points on $[-1, 1]$ as

$$t_n = \frac{t_f}{2}\left(1 - \cos\left(\frac{\pi n}{N}\right)\right), \quad n = 0, 1, 2, \ldots, N. \tag{7}$$

The discrete orthogonality relation is vital for the ChSCM defined as

$$\sum_{n=0}^{N} \alpha_n \bar{T}_i(t_n) \bar{T}_j(t_n) = \gamma_i \delta_{ij} \tag{8}$$

with

$$\gamma_i = \begin{cases} \frac{N}{2}, & i \neq 0, N \\ N, & i = 0, N. \end{cases} \tag{9}$$

With the use of discrete orthogonality relations (7), we get the following equality

$$c_j(x, y) = \frac{2}{N} \sum_{n=0}^{N} \alpha_n \bar{T}_j(t_n) u(x, y, t_n). \tag{10}$$

Thus, we can calculate the $u_t(x, y, t)$ at the predetermined collocation points as follows,

$$u_t(x, y, t_i) = \sum_{n=0}^{N} [A_t]_{in} c_n(x, y), \tag{11}$$

where $i = 0, 1, 2, \ldots, N$, $c_n(x, y) = u(x, y, t_n)$ and $A_t$ is $(N \times N)$ matrix defined as

$$[A_t]_{in} = \begin{cases} \mu \sum_{j=0}^{N} (-1)^j j^2 \cos\left(j\left(\pi - \frac{(N+n-1)\pi}{N}\right)\right), & i = 0 \\ \mu \sum_{j=0}^{N} \frac{j \sin\left(j\left(\pi - \frac{(i-1)\pi}{N}\right)\right) \cos\left(j\left(\pi - \frac{(n-1)\pi}{N}\right)\right)}{\sqrt{1-\cos^2\left(\pi - \frac{(i-1)\pi}{N}\right)}}, & i = 1, 2, \ldots, N-1 \\ \mu \sum_{j=0}^{N} (-1)^{j+1} j^2 \cos\left(j\left(\pi - \frac{(n-1)\pi}{N}\right)\right), & i = N \end{cases} \tag{12}$$

where $\mu = \frac{2}{t_f}$ and $[A_t]_{in}$ is called as the Chebyshev differentiation matrix. Inserting the approximation function (5) into main equation (1) in time collocation sense, we obtain the following system of elliptic boundary value problems

$$D_x \frac{\partial^2 c}{\partial x^2} + D_y \frac{\partial^2 c}{\partial y^2} - V_x \frac{\partial c}{\partial x} - V_y \frac{\partial c}{\partial y} - A_t c = 0 \qquad (13)$$

where $c(x,y) = [c_0(x,y), c_1(x,y), c_2(x,y), \ldots, c_N(x,y)]^T$ is $((N+1) \times 1)$ column vector. Applying the initial condition $u(x,y,0) = c_0(x,y) = g(x,y)$ reduces the system into the following form

$$D_x \frac{\partial^2 c}{\partial x^2} + D_y \frac{\partial^2 c}{\partial y^2} - V_x \frac{\partial c}{\partial x} - V_y \frac{\partial c}{\partial y} - \bar{A}_t c = G(x,y) \qquad (14)$$

where $G(x,y)$ occurs due to the initial condition (3) and $\bar{A}_t$ is $(N \times N)$ matrix reduced from $A_t$. Thus, equation (14) represents the $(N \times N)$ system of boundary value problems with the following boundary conditions,

$$c(a,y) = f_1(y,\bar{t}) \text{ and } c(b,y) = f_2(y,\bar{t}) \qquad (15)$$

$$c(x,c) = h_1(x,\bar{t}) \text{ and } c(x,d) = h_2(x,\bar{t}) \qquad (16)$$

where $\bar{t} = [t_1, t_2, \ldots, t_N]^T$ defined in equation (6).

## 3.2 The IELDTM for Two Dimensions

In this subsection, the IELDTM will be derived to solve 2D system of boundary value problems (14). The detailed background of the differential transformation for multi-dimensional analytic functions can be seen in literature [21]. For simplicity, we omit these details here.

Let us divide the spatial domain $A = [a,b] \times [c,d]$ into $M_x M_y$ spatial element $[x_i, x_{i+1}] \times [y_j, y_{j+1}]$ with $\Delta x = x_{i+1} - x_i$, $\Delta y = y_{j+1} - y_j$ and the centre of each element is defined at $\left(x_{i+\frac{1}{2}}, y_{j+\frac{1}{2}}\right)$. The function $c(x,y)$ can be locally represented with the convergent local Taylor expansion about $\left(x_{i+\frac{1}{2}}, y_{j+\frac{1}{2}}\right)$ as follows

$$c_{ij}(x,y) = \sum_{k=0}^{K} \sum_{p=0}^{K-k} C_{ij}(k,p) \left(x - x_{i+\frac{1}{2}}\right)^k \left(y - y_{j+\frac{1}{2}}\right)^p + O\left(\left(x - x_{i+\frac{1}{2}}\right)^{K+1} + \left(y - y_{j+\frac{1}{2}}\right)^{K+1}\right) \qquad (17)$$

where $C_{ij}(k,p) = \frac{1}{k!p!} \frac{\partial^{k+p} c(x,y)}{\partial x^k \partial y^p}\big|_{x=x_{i+\frac{1}{2}}, y=y_{j+\frac{1}{2}}}$ is the local differential transform of the function $c(x,y)$ about $x = x_{i+\frac{1}{2}}$ and $y = y_{j+\frac{1}{2}}$ for $i = 0,1,\ldots,M_x - 1$ and $j = 0,1,\ldots,M_y - 1$. To completely describe the convergent local solution $c_{ij}(x,y)$ of order $K$, equation (17) includes

$N\frac{(K+1)(K+2)}{2}$ coefficients to be determined. The crucial advantage of the differential transformation leads to the reduction of these unknown coefficients to $N(2K+1)$. Taking the differential transform of equation (14) yields

$$C_{ij}(k,p+2) = \frac{1}{D_y(p+1)(p+2)}\left[-D_x(k+1)(k+2)C_{ij}(k+2,p) + V_x(k+1)C_{ij}(k+1,p) + V_y(p+1)C_{ij}(k,p+1) - \bar{A}_t C_{ij}(k,p) + G_{ij}(k,p)\right] \quad (18)$$

where , $p = 0,1,\ldots,K-1$, $i = 0,1,\ldots,M_x-1$, $j = 0,1,\ldots,M_y-1$ and $G_{ij}(k,p)$ is the local differential transform of the $G(x,y)$ at $x = x_{i+\frac{1}{2}}$ and $y = y_{j+\frac{1}{2}}$. Thus, the use of recursion (18), one only needs to determine the following set of $N(2K+1)$ local values

$$V_{ij} = \{C_{ij}(k,0)|\ k \leq K \wedge k \in N\} \cup \{C_{ij}(k,1)|\ k \leq K-1 \wedge k \in N\} \quad (19)$$

where $i = 0,1,\ldots,M_x - 1$ and $j = 0,1,\ldots,M_y - 1$. Thus, the total degree of freedom is $NM_xM_y(2K+1)$. As proven in numerical experiments, the minimization of the local degrees of freedom leads to higher-order approximations with optimum costs.

The function $c_{ij}(x,y)$ is assumed as analytic for all $(i,j)$ and has the radius of convergences $\rho_x^i > \Delta x$ and $\rho_y^j > \Delta y$. This assumption leads us to search for the relations between neighbour solutions $c_{ij}(x,y)$, $c_{i(j+1)}(x,y)$ and $c_{(i+1)j}(x,y)$ in terms of their local behaviours. Let us define the intervals $P^i = [x_{i+\frac{1}{2}} - \rho_x^i, x_{i+\frac{1}{2}} + \rho_x^i]$, $P^{i+1} = [x_{i+\frac{3}{2}} - \rho_x^{i+1}, x_{i+\frac{3}{2}} + \rho_x^{i+1}]$, $W^j = [y_{j+\frac{1}{2}} - \rho_y^j, y_{j+\frac{1}{2}} + \rho_y^j]$, $W^{j+1} = [y_{j+\frac{3}{2}} - \rho_y^{j+1}, y_{j+\frac{3}{2}} + \rho_y^{j+1}]$, $P = [x_{i+\frac{1}{2}}, x_{i+\frac{3}{2}}]$. and $W = [y_{j+\frac{1}{2}}, y_{j+\frac{3}{2}}]$. Assume that the conditions $P \subset (P^i \cap P^{i+1})$ and $W \subset (W^j \cap W^{j+1})$ are hold, then we state the following conclusions;

- $C^0$ and $C^1$ continuity conditions can be applied over the region $[x_{i+\frac{1}{2}}, x_{i+\frac{3}{2}}] \times [y_{j+\frac{1}{2}}, y_{j+\frac{3}{2}}]$.
- Any point on the interval $[x_{i+\frac{1}{2}}, x_{i+\frac{3}{2}}]$ can be written as $x^* = x_{i+\frac{1}{2}} + (1-\theta_x)\Delta x$, where $0 \leq \theta_x \leq 1$.
- Any point on the interval $[y_{j+\frac{1}{2}}, x_{j+\frac{3}{2}}]$ can be written as $y^* = y_{j+\frac{1}{2}} + (1-\theta_y)\Delta y$, where $0 \leq \theta_y \leq 1$.

We define the sets $H_i^x = \{x_i, x_i + dx, \ldots, x_{i+1}\}$ and $H_j^y = \{y_j, y_j + dy, \ldots, y_{j+1}\}$ with $dx = \frac{\Delta x}{S}$ and $dy = \frac{\Delta y}{S}$ where $S$ denotes the partition number of each edge. The sets $H_i^x$ and $H_j^y$ include all possible collocation points in $x$ and $y$ directions, respectively. $C^0$ and $C^1$ continuity equations in $x - y$ directions and boundaries yield

$$c_{ij}(x,y) = c_{(i+1)j}(x,y) \text{ and } \frac{\partial c_{ij}(x,y)}{\partial n} = \frac{\partial c_{(i+1)j}(x,y)}{\partial n} \text{ for } \forall (x,y) \in W_x^{ij}, \quad (20)$$

$$c_{ij}(x,y) = c_{i(j+1)}(x,y) \text{ and } \frac{\partial c_{ij}(x,y)}{\partial n} = \frac{\partial c_{i(j+1)}(x,y)}{\partial n} \text{ for } \forall (x,y) \in W_y^{ij}, \quad (21)$$

$$c_{ij}(x,y) = F_{ij}(x,y) \text{ and } \frac{\partial c_{ij}(x,y)}{\partial n} = \frac{\partial F_{ij}(x,y)}{\partial n} \text{ for } \forall (x,y) \in \partial A_{ij}, \quad (22)$$

where $W_x^{ij} = \left\{(x,y) | \ x = x_{i+\frac{1}{2}} + (1-\theta_x)\Delta x \wedge y \in H_j^y\right\}$, $W_y^{ij} = \{(x,y) | \ x \in H_i^x \wedge y = y_{j+\frac{1}{2}} + (1-\theta_y)\Delta y\}$, $\partial A_{ij}$ denotes the discretized computational boundary of the domain $A$ and the $F_{ij}(x,y)$ is the discrete boundary data given in equations (15)-(16). Thus, combining all continuity equations defined in equations (20)-(22) leads to the following linear algebraic system

$$M\zeta = P \quad (23)$$

where $M$ is $\left(N[6M_xM_y + M_x + M_y](S+1)\right) \times \left(NM_xM_y(2K+1)\right)$ matrix, $P$ is $\left(N[6M_xM_y + M_x + M_y](S+1)\right) \times 1$ column vector and $\zeta$ is $\left(NM_xM_y(2K+1)\right) \times 1$ column vector defined as

$$\zeta = \left[V_{00}, V_{01}, \dots, V_{(M_x-1)(M_y-1)}\right]^T \quad (24)$$

where $V_{ij}$ is $1 \times (N(2K+1))$ vector defined in equation (19). By suitable selection of $S$ with $S > \frac{M_xM_y(2K+1)}{6M_xM_y+M_x+M_y}$, equation (23) describes an overdetermined linear system. By minimizing the residual norm $\|M\zeta - P\|$ in the least square sense, all unknown values $V_{ij}$ are obtained. Note that $\theta_x$ and $\theta_x$ determines the direction of the continuity equations, i.e. the schemes are:

- forward directional with $\theta_x = 0$ and $\theta_y = 0$
- backward directional with $\theta_x = 1$ and $\theta_y = 1$
- central directional with $\theta_x = 1/2$ and $\theta_y = 1/2$.

As proven in the later sections, the central directional schemes are more accurate than the other selections.

## 3. IELDTM-ChSCM for (3+1) dimensional AD equation

In this section, we introduce the IELDTM-ChSCM hybridization for solving the following (3+1) dimensional advection-diffusion (AD) equation

$$u_t + V_x u_x + V_y u_y + V_z u_z = D_x u_{xx} + D_y u_{yy} + D_z u_{zz}, \quad (x,y,z) \in A \text{ and } t \in [0, t_f] \quad (25)$$

with Dirichlet boundary conditions

$$u(a,y,z,t) = f_1(y,z,t) \text{ and } u(b,y,z,t) = f_2(y,z,t) \tag{26}$$

$$u(x,c,z,t) = h_1(x,z,t) \text{ and } u(x,d,z,t) = h_2(x,z,t) \tag{27}$$

$$u(x,y,e,t) = w_1(x,y,t) \text{ and } u(x,y,f,t) = w_2(x,y,t) \tag{27}$$

and initial condition

$$u(x,y,z,0) = g(x,y,z) \tag{28}$$

where $A = [a,b] \times [c,d] \times [e,f]$, $V_x$, $V_y$ and $V_z$ are advection constants, $D_x, D_y$ and $D_z$ are diffusion constants and $f_1, f_2, h_1, h_2, w_1, w_2$ and $g$ are known smooth functions. The subscripts $x, y, z,$ and $t$ represent differentiations with respect to spaces and time, respectively. As we constructed for the (2+1) dimensional AD equation (1), (3+1) dimensional AD equation (25) will be reduced to a system of three-dimensional elliptic boundary value problems with the use of the Chebyshev spectral collocation method (ChSCM) in time.

### 3.1. Temporal variation for (3+1) dimensional AD equation

Since all the details of the ChSCM have been given in Section 2.1, just the essential formulae are presented here. Let us approximate the time part of $u(x,y,z,t)$ with $Nth$ order Chebyshev polynomials as

$$u(x,y,z,t) = \sum_{n=0}^{N} \alpha_n c_n(x,y,z) \bar{T}_n(t) \tag{29}$$

where $c_n(x,y,z)$ represent the spatial parts of the $u(x,y,z,t)$. The use of approximation function (29) into (3+1) dimensional AD equation (25) leads to the following three-dimensional system of elliptic boundary value problems

$$D_x \frac{\partial^2 \mathbf{c}}{\partial x^2} + D_y \frac{\partial^2 \mathbf{c}}{\partial y^2} + D_y \frac{\partial^2 \mathbf{c}}{\partial z^2} - V_x \frac{\partial \mathbf{c}}{\partial x} - V_y \frac{\partial \mathbf{c}}{\partial y} - V_z \frac{\partial \mathbf{c}}{\partial z} - A_t \mathbf{c} = \mathbf{0}$$

where $\mathbf{c}(x,y,z) = [c_0(x,y,z), c_1(x,y,z), c_2(x,y,z), \ldots, c_N(x,y,z)]^T$ is $((N+1) \times 1)$ column vector. Applying the initial condition $u(x,y,z,0) = c_0(x,y,z) = g(x,y,z)$ reduces the system into the following form

$$D_x \frac{\partial^2 \mathbf{c}}{\partial x^2} + D_y \frac{\partial^2 \mathbf{c}}{\partial y^2} + D_z \frac{\partial^2 \mathbf{c}}{\partial z^2} - V_x \frac{\partial \mathbf{c}}{\partial x} - V_y \frac{\partial \mathbf{c}}{\partial y} - V_z \frac{\partial \mathbf{c}}{\partial z} - \bar{A}_t \mathbf{c} = \mathbf{G}(x,y,z) \tag{30}$$

where $\mathbf{G}(x,y,z)$ occurs due to the initial condition (28) and $\bar{A}_t$ is $(N \times N)$ matrix reduced from $A_t$. Thus, equation (30) represents the $(N \times N)$ system of nonhomogeneous three-dimensional boundary value problems with the following boundary conditions,

$$\mathbf{c}(a,y,z) = f_1(y,z,\bar{t}) \text{ and } \mathbf{c}(b,y,z) = f_2(y,z,\bar{t}) \tag{31}$$

$$\mathbf{c}(x,c,z) = h_1(x,z,\bar{t}) \text{ and } \mathbf{c}(x,d,z) = h_2(x,z,\bar{t}) \tag{32}$$

$$\mathbf{c}(x,y,e) = w_1(x,y,\bar{t}) \text{ and } \mathbf{c}(x,y,f) = w_2(x,y,\bar{t}) \tag{33}$$

where $\bar{t} = [t_1, t_2, \ldots, t_N]^T$ defined in equation (6).

## 3.2 The IELDTM for Three Dimensions

The IELDTM is derived here to solve 3D system of boundary value problems (30). In Section 2.2, the method was explained with all details, and here we omit these details for simplicity. Let us divide the spatial domain $A = [a, b] \times [c, d] \times [e, f]$ into $M_x M_y M_z$ spatial element $[x_i, x_{i+1}] \times [y_j, y_{j+1}] \times [z_r, z_{r+1}]$ with $\Delta x = x_{i+1} - x_i$, $\Delta y = y_{j+1} - y_j$, $\Delta z = z_{r+1} - z_r$ and the centre of each element is defined at $\left(x_{i+\frac{1}{2}}, y_{j+\frac{1}{2}}, z_{r+\frac{1}{2}}\right)$. The function $c(x, y, z)$ can be locally represented with the convergent local Taylor expansion about $\left(x_{i+\frac{1}{2}}, y_{j+\frac{1}{2}}, z_{r+\frac{1}{2}}\right)$ as follows

$$c_{ijr}(x,y,z) = \sum_{h=0}^{K} \sum_{p=0}^{h} \sum_{k=0}^{h-p} C_{ijr}(k, p, h-p-k) \left(x - x_{i+\frac{1}{2}}\right)^k \left(y - y_{j+\frac{1}{2}}\right)^p \left(z - z_{r+\frac{1}{2}}\right)^{h-p-k} + O\left(\left(x - x_{i+\frac{1}{2}}\right)^{K+1} + \left(y - y_{j+\frac{1}{2}}\right)^{K+1} + \left(z - y_{r+\frac{1}{2}}\right)^{K+1}\right) \quad (34)$$

where $C_{ijr}(k, p, h) = \frac{1}{k!p!h!} \frac{\partial^{k+p+h} c(x,y,z)}{\partial x^k \partial y^p \partial z^h}\Big|_{x=x_{i+\frac{1}{2}}, y=y_{j+\frac{1}{2}}, z=z_{r+\frac{1}{2}}}$ is the local differential transform of the function $c(x, y, z)$ about $x = x_{i+\frac{1}{2}}, y = y_{j+\frac{1}{2}}$ and $z = z_{r+\frac{1}{2}}$ for $i = 0, 1, \ldots, M_x - 1$, $j = 0, 1, \ldots, M_y - 1$ and $r = 0, 1, \ldots, M_z - 1$. The convergent local solution $c_{ijr}(x, y, z)$ of order $K$ includes $N\frac{(K+1)(K+2)(K+3)}{6}$ coefficients to be determined. With the use of the differential transformation of three-dimensional equation (30), the number of unknowns will be reduced to $N(K + 1)^2$. Taking the differential transform of equation (30) yields

$$C_{ijr}(k, p, h+2) = \frac{1}{D_z(h+1)(h+2)} \big[-D_x(k+1)(k+2)C_{ijr}(k+2, p, h) - D_y(p+1)(p+2)C_{ijr}(k, p+2, h) + V_x(k+1)C_{ijr}(k+1, p, h) + V_y(p+1)C_{ijr}(k, p+1, h) + V_z(h+1)C_{ijr}(k, p, h+1) - \bar{A}_t C_{ijr}(k, p, h) + G_{ijr}(k, p, h)\big] \quad (35)$$

where $k, p, h = 0, 1, \ldots, K - 1$ and $G_{ijr}(k, p, r)$ is the local differential transform of the $G(x, y, z)$. Thus, with the use of recursion (35), one only needs to determine the following set of $N(K + 1)^2$ local values

$$V_{ijr} = \{C_{ijr}(k, p, 0) | \; k + p \le K \wedge k, p \in \mathbb{N}\} \cup \{C_{ijr}(k, p, 1) | \; k + p \le K - 1 \wedge k, p \in \mathbb{N}\} \quad (36)$$

where $i = 0, 1, \ldots, M_x - 1$, $j = 0, 1, \ldots, M_y - 1$ and $r = 0, 1, \ldots, M_z - 1$. The total degree of freedom is $M_x M_y M_z N(K + 1)^2$.

We define the set of collocation points $H_i^x = \{x_i, x_i + dx, \ldots, x_{i+1}\}$, $H_j^y = \{y_j, y_j + dy, \ldots, y_{j+1}\}$ and $H_r^z = \{z_r, z_r + dz, \ldots, z_{r+1}\}$ with $dx = \frac{\Delta x}{S}$, $dy = \frac{\Delta y}{S}$ and $dz = \frac{\Delta z}{S}$ where $S$ denotes the partition number of each edge. $C^0$ and $C^1$ continuity equations in $x$, $y$, and $z$ directions and boundaries lead to

$$\boldsymbol{c}_{ijr}(x,y,z) = \boldsymbol{c}_{(i+1)jr}(x,y,z) \text{ and } \frac{\partial c_{ijr}(x,y,z)}{\partial n} = \frac{\partial c_{(i+1)jr}(x,y,z)}{\partial n} \text{ for } \forall (x,y,z) \in W_x^{ijr}, \quad (37)$$

$$\boldsymbol{c}_{ijr}(x,y,z) = \boldsymbol{c}_{i(j+1)r}(x,y,z) \text{ and } \frac{\partial c_{ijr}(x,y,z)}{\partial n} = \frac{\partial c_{i(j+1)r}(x,y,z)}{\partial n} \text{ for } \forall (x,y,z) \in W_y^{ijr}, \quad (38)$$

$$\boldsymbol{c}_{ijr}(x,y,z) = \boldsymbol{c}_{ij(r+1)}(x,y,z) \text{ and } \frac{\partial c_{ijr}(x,y,z)}{\partial n} = \frac{\partial c_{ij(r+1)}(x,y,z)}{\partial n} \text{ for } \forall (x,y,z) \in W_z^{ijr}, \quad (39)$$

$$\boldsymbol{c}_{ijr}(x,y,z) = \boldsymbol{F}_{ijr}(x,y,z) \text{ and } \frac{\partial c_{ijr}(x,y,z)}{\partial n} = \frac{\partial F_{ijr}(x,y,z)}{\partial n} \text{ for } \forall (x,y,z) \in \partial A_{ijr}, \quad (40)$$

where $W_x^{ijr} = \left\{(x,y,z) | x = x_{i+\frac{1}{2}} + (1-\theta_x)\Delta x \wedge y \in H_j^y \wedge z \in H_r^z\right\}$, $W_y^{ijr} = \Big\{(x,y,z) | x \in H_i^x \wedge y = y_{j+\frac{1}{2}} + (1-\theta_y)\Delta y \wedge z \in H_r^z\Big\}$, $W_z^{ijr} = \Big\{(x,y,z) | x \in H_i^x \wedge y \in H_j^y \wedge z = z_{r+\frac{1}{2}} + (1-\theta_z)\Delta z\Big\}$, $\partial A_{ijr}$ represents the discretized computational boundary of the domain $A$ and $\boldsymbol{F}_{ijr}(x,y,z)$ is the discrete boundary data that are given in equations (31)-(33). The direction parameters $0 \leq \theta_x, \theta_y, \theta_z \leq 1$ determine the implicit- explicit continuity relations and have a crucial role in the accuracy of the numerical algorithm, as discussed in numerical experiments. Assembling all continuity equations (37)-(40) yields to the following linear algebraic system

$$M\boldsymbol{\zeta} = P \quad (41)$$

where $M$ is $\left(N[9M_xM_yM_z + M_xM_y + M_xM_z + M_yM_z](S+1)^2\right) \times \left(NM_xM_yM_z(K+1)^2\right)$ matrix, $P$ is $\left(N[9M_xM_yM_z + M_xM_y + M_xM_z + M_yM_z](S+1)^2\right) \times 1$ column vector and $\boldsymbol{\zeta}$ is $\left(NM_xM_yM_z(K+1)^2\right) \times 1$ column vector defined as

$$\boldsymbol{\zeta} = \left[V_{000}, V_{001}, \dots, V_{(M_x-1)(M_y-1)(M_z-1)}\right]^T \quad (42)$$

where $V_{ij}$ is $1 \times (N(K+1)^2)$ vector defined in equation (36). If $S$ is selected as $(S+1)^2 > \frac{(NM_xM_yM_z(K+1)^2)}{9M_xM_yM_z+M_xM_y+M_xM_z+M_yM_z}$, equation (43) describes an overdetermined linear system and minimizing the residual norm $\|M\boldsymbol{\zeta} - P\|$ in the least square sense will lead to unknown values $V_{ij}$. Note that, $\theta_x$, $\theta_y$ and $\theta_z$ are the direction parameters in $x$, $y$, and $z$ directions, respectively.

## 4 Numerical Experiments

In this section, the IELDTM derived in Sections 2-3 is tested over the linear/nonlinear (2+1) dimensional AD problems and a linear (3+1) dimensional AD problem. The convergence of the present algorithms is quantitatively demonstrated with the consideration of both $h$ −refinement and $p$ −refinement procedures. The produced numerical results are compared with the finite difference, finite element, and spectral presented produced in literature [2,8,22,23]. To evaluate error norms of the present results, we prefer to use the following definitions,

$$E_{abs}^{ij} = |u_{ij}^{exact} - u_{ij}^{numerical}|,$$

$$\|E\|_\infty = \max_{(x,y)\in D, t\in[0,t_f]} |u^{exact}(x,y,t) - u^{numerical}(x,y,t)|.$$

**Problem 1 [8,23]**

Consider the (2+1) dimensional pure diffusion process with the choices of $V_x = V_y = 0$, $D_x = D_y = 1$ in the AD equation (1) for which the exact solution is given by [8]

$$u(x,y,t) = e^{-2\pi^2 t}\sin(\pi x)\sin(\pi y), t > 0 \text{ and } (x,y) \in [0,1] \times [0,1]. \tag{43}$$

The homogeneous Dirichlet boundary conditions and the initial condition are taken from the exact solution (34) as follows:

$$u(0,y,t) = u(1,y,t) = 0, \tag{44}$$

$$u(x,0,t) = u(x,1,t) = 0, \tag{45}$$

$$u(x,y,0) = \sin(\pi x)\sin(\pi y). \tag{46}$$

In Table 1, the present IELDTM-ChSCM algorithm has been compared with the finite difference-based methods [8,23] with the consideration of maximum error norms. The spatial degrees of freedom (*dof*) of all numerical algorithms have been presented comparatively. As seen in Table 1, the IELDTM with forward, backward, and central directions produce more accurate results with far less spatial *dof* than the finite difference-based methods presented in [8, 23]. In Figure 1, pointwise absolute errors produced by the central IELDTM are presented at various times with contour plots. As clearly observed from the figure, the IELDTM provides highly accurate results with optimized degrees of freedom. Even if we consider only N = 20 elements in time, the ChSCM offers continuous time integration of the advection-diffusion equation up to $t = 1$ by maintaining roughly pointwise absolute errors about $10^{-9}$.

**Problem 2 [8]**

Consider the (2+1) dimensional advection-diffusion process with $V_x = V_y = 1$ for which the exact solution is given by [8]

$$u(x,y,t) = \frac{1}{4t+1}\exp\left(-\frac{(x-t-0.5)^2}{D_x(4t+1)} - \frac{(y-t-0.5)^2}{D_y(4t+1)}\right), t > 0 \text{ and } (x,y) \in [0,x_f] \times [0,y_f]. \tag{47}$$

The nonhomogeneous time-dependent Dirichlet boundary conditions can be taken from exact solution (47), and the initial condition is considered as

$$u(x,y,0) = exp\left(-\frac{(x-0.5)^2}{D_x} - \frac{(y-0.5)^2}{D_y}\right), (x,y) \in [0,x_f] \times [0,y_f]. \tag{48}$$

**Table 1** Comparison of the maximum error norms $\|E\|_\infty$ of various present schemes with the FDM-based results produced in [8, 23] at $t_f = 0.25$, $N = 15$, $S = 14$ and $K = 10$ for Problem 1.

| FDM | | | IELDTM | | | |
|---|---|---|---|---|---|---|
| $M_x \times M_y$ (dof) | [23] | [8] | $M_x \times M_y$ (dof) | $\theta_x = \theta_y = 0.5$ | $\theta_x = \theta_y = 1$ | $\theta_x = \theta_y = 0$ |
| 11 × 11 (121) | 1.36E-05 | 5.67E-06 | 2 × 2 (84) | 4.14E-08 | 6.91E-06 | 6.91E-06 |
| 21 × 21 (441) | 8.55E-07 | 3.36E-07 | 3 × 3 (189) | 1.69E-09 | 1.68E-06 | 1.68E-06 |
| 41 × 41 (1681) | 5.34E-08 | 2.07E-08 | 4 × 4 (336) | 8.70E-11 | 3.63E-08 | 3.63E-08 |
| 81 × 81 (6561) | 3.34E-09 | 1.29E-09 | 5 × 5 (525) | 9.99E-12 | 1.81E-09 | 1.81E-09 |

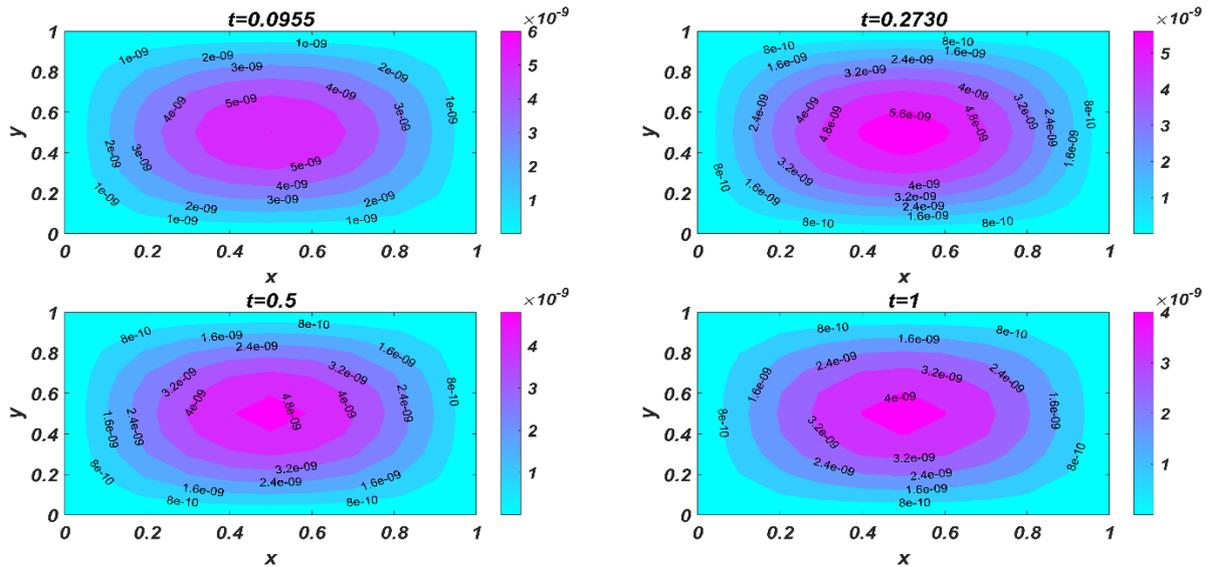

**Figure 1:** Descriptive behavior of pointwise absolute errors produced by the central IELDTM at various times with $\theta_x = \theta_y = 0.5$, $M_x = M_y = 5$, $K = 10$, $N = 20$, $S = 14$ and $t_f = 1$.

In Figures 2-3, spectral convergence of the IELDTM with has been proven with respect to $K-$refinement (order refinement) and $\Delta x/\Delta y-$refinement for various values of the direction parameters $\theta_x$ and $\theta_y$. As theoretically expected, $K-$refinement procedure leads to an exponential convergence and $\Delta x/\Delta y-$refinement procedure leads to a polynomial convergence like $\Delta h^P$ when $\Delta x = \Delta y = \Delta h$. Since $K-$refinement procedure increases the local degrees of freedom (*ldof*) linearly as $ldof = 2K + 1$, the IELDTM provides highly accurate results with optimized *dof*. The mesh discretization parameter $S$ plays a vital role in both the accuracy and the computational cost of the IELDTM. In Figure 4a, the maximum error norms produced by the IELDTM have been presented with varying values of $S$. A rapid decrease in the corresponding error norm has been observed, and the curve behaves asymptotically for higher values of $S$. The mesh discretization ratio $\phi = Number\ of\ equations/Number\ of\ unknowns$ versus the maximum error norm curve has been illustrated with the varying S values in Figure 4b. The curve indicates that the errors begin to be stable when the matrix that occurred in equation (23) is almost square ($\phi \cong 1$). Advection-dominated fluid flow problems are extremely challenging for all numerical methods [1]. With the consideration of various diffusion coefficients $(D_x, D_x)$ in both space dimensions, the central IELDTM solutions and the corresponding pointwise absolute errors are illustrated in Figure 5. By considering the parameter values for $\theta_x = \theta_y = 0.5$, $x_f = y_f = 1$, $K = 4$, $M_x = M_y = 18$, $N = 5$, $S = 6$ and $t_f = 0.05$, it can be observed that the present algorithm produces accurate results for both non-stiff and stiff cases. The performance of the present IELDTM is illustrated over a relatively large space-time computational domain with $x_f = y_f = 5$ and $t_f = 2$ in Figure 6. As observed from the figure, the IELDTM accurately captures the physical behaviour with optimized degrees of freedom.

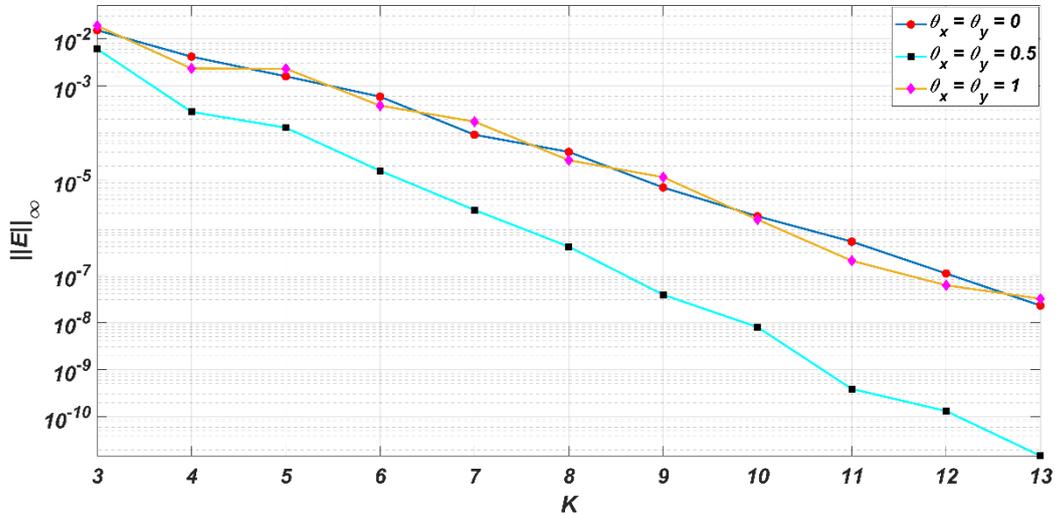

**Figure 2:** Convergence results of the present algorithm with respect to spectral convergence order $K$ for various values of the direction parameter couples $(\theta_x, \theta_y)$, $x_f = y_f = 1$, $D_x = D_y = 1$, $M_x = M_y = 2$, $N = 10$, $S = 10$ and $t_f = 0.1$.

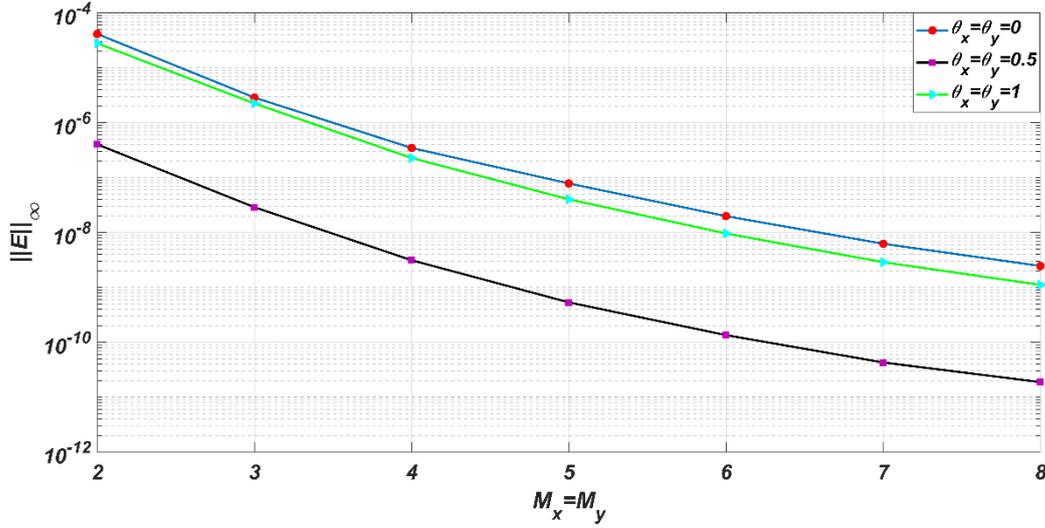

**Figure 3:** Convergence results of the present algorithm with respect to the spatial element numbers $M_x = M_y$ for various values of the direction parameter couples $(\theta_x, \theta_y)$, $x_f = y_f = 1$, $D_x = D_y = 1$, $K = 8$, $N = 10$, $S = 12$ and $t_f = 0.1$.

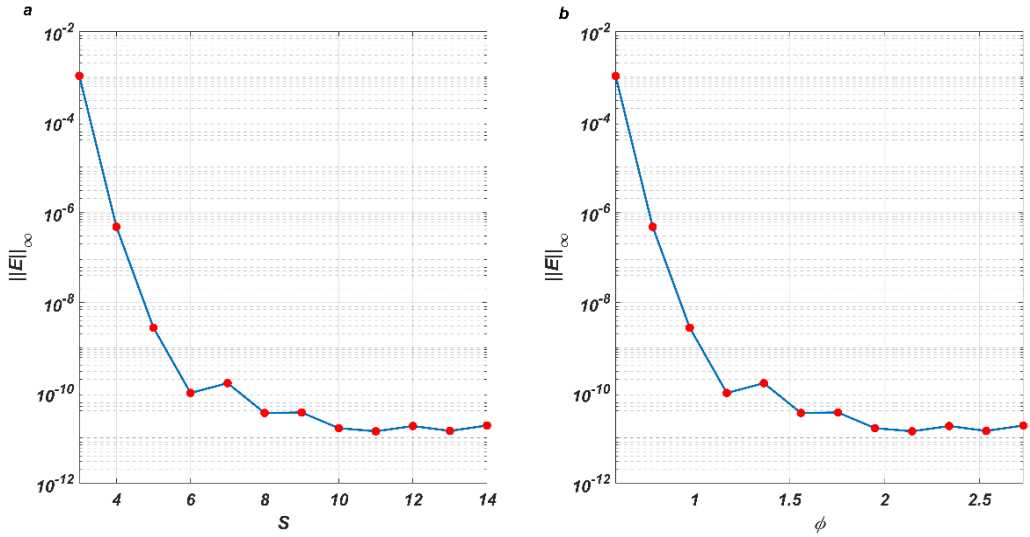

**Figure 4:** a) Effects of the mesh discretization parameter $S$ to the maximum error norms $\|E\|_\infty$, b) Effects of the mesh discretization ratio $\phi$ to the maximum error norms $\|E\|_\infty$ with $x_f = y_f = 1$, $D_x = D_y = 1$, $M_x = M_y = 1$, $N = 10$, $K = 20$ and $t_f = 0.1$.

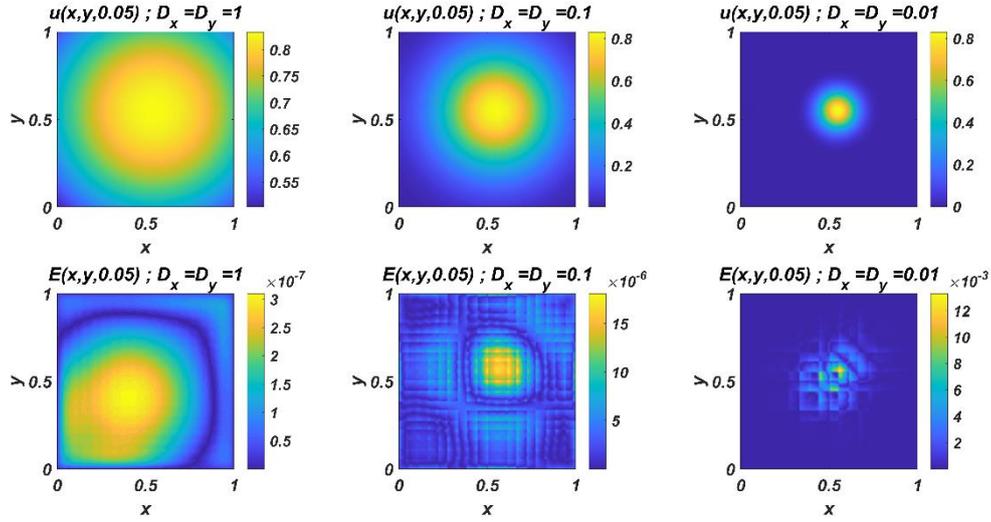

**Figure 5:** Numerical results and corresponding pointwise absolute errors for solving Problem 2 with the use of the central IELDTM for $\theta_x = \theta_y = 0.5$, $x_f = y_f = 1$, $K = 4$, $M_x = M_y = 18$, $N = 5$, $S = 6$ and $t_f = 0.05$.

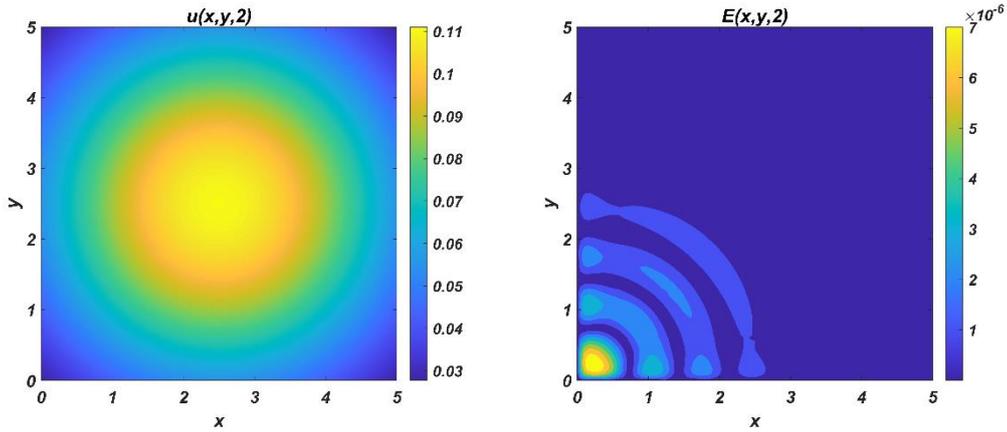

**Figure 6:** Numerical results and corresponding pointwise absolute errors for solving Problem 2 with the use of the central IELDTM for $\theta_x = \theta_y = 0.5$, $x_f = y_f = 5$, $K = 8$, $M_x = M_y = 8$, $N = 16$, $S = 5$ and $t_f = 2$.

**Problem 3 [2,22]**

Let us consider the following (2+1) dimensional nonlinear Burgers equation

$$u_t + uu_x + uu_y = D(u_{xx} + u_{yy}), \quad (x,y) \in [a,b] \times [c,d], \quad t \in [0, t_f] \tag{49}$$

with the following exact solution [2]

$$u(x,y,t) = \frac{1}{1+e^{(x+y-t)/D}} \tag{50}$$

where $D$ is kinematic viscosity constant along with the $x$ and $y$ directions. The initial and Dirichlet boundary conditions can be taken from the exact solution (50).

In Table 2, the present IELDTM with backward, central, and forward directions has been compared with the finite element method [2] and the Chebyshev spectral collocation method [22] by measuring the maximum error norms for solving (2+1) dimensional nonlinear advection-diffusion equation (49). As seen in Table 2, the proposed method offers much more accurate results using far less *dof* than both techniques. The spectral convergence results for solving nonlinear equation (49) with respect to K-refinement and $\Delta x/\Delta y$ −refinement procedures have been illustrated in Figures 7-8. The expected convergence behaviors are seen to be satisfied.

**Table 2** Comparison of the maximum error norms $\|E\|_\infty$ of various present schemes with the FEM and the ChSCM produced in [2,22] at $t_f = 0.25, N = 15, S = 16, D = 1$ and $K = 10$ for Problem 3.

| Literature | | | IELDTM | | | |
| --- | --- | --- | --- | --- | --- | --- |
| $M_x \times M_y$ (*dof*) | ChSCM [22] | FEM [2] | $M_x \times M_y$ (*dof*) | $\theta_x = \theta_y = 0.5$ | $\theta_x = \theta_y = 1$ | $\theta_x = \theta_y = 0$ |
| $5 \times 5$ (25) | 8.94E-08 | 4.65E-08 | $1 \times 1$ (21) | 5.84E-10 | 5.84E-10 | 5.84E-10 |
| $10 \times 10$ (100) | 7.45E-07 | 5.90E-09 | $2 \times 2$ (84) | 1.12E-14 | 1.11E-12 | 6.46E-13 |
| $15 \times 15$ (225) | - | 2.18E-09 | $3 \times 3$ (189) | 5.55E-17 | 4.18E-14 | 4.08E-14 |
| $30 \times 30$ (900) | - | 1.01E-09 | $4 \times 4$ (336) | 1.66E-16 | 3.83E-15 | 4.21E-15 |

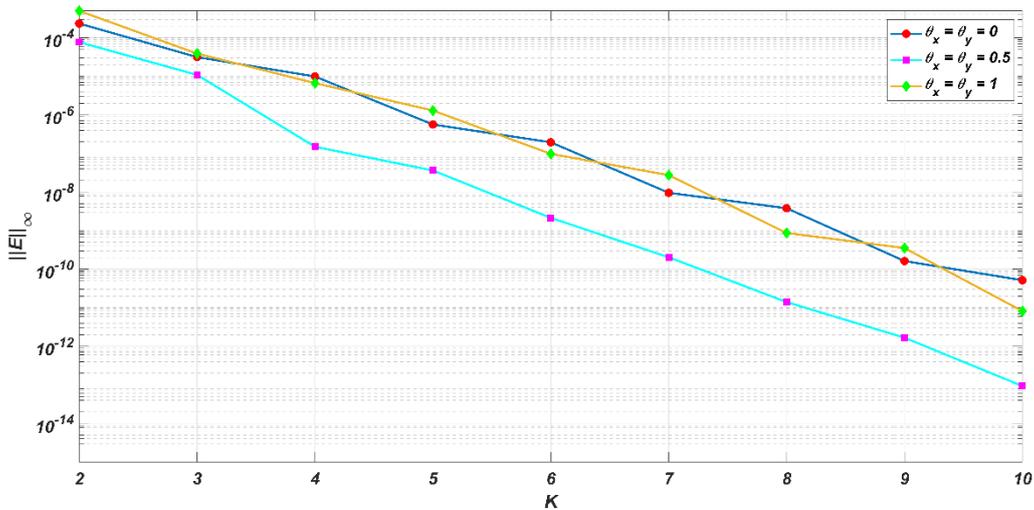

**Figure 7:** Convergence results of the present algorithm with respect to spectral convergence order $K$ for various values of the direction parameter couples $(\theta_x, \theta_y)$, $M_x = M_y = 2, N = 10, S = 10, D = 1$ and $t_f = 0.5$.

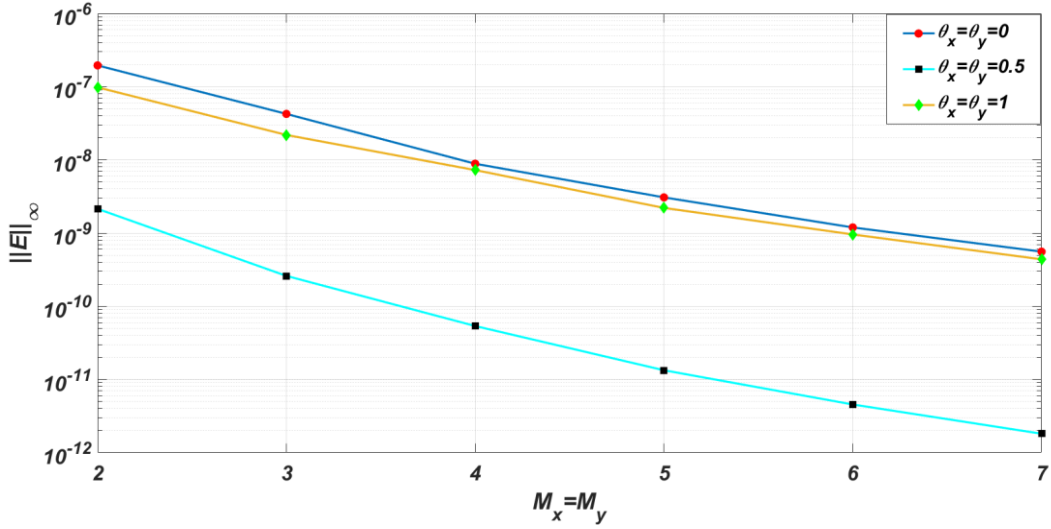

**Figure 8:** Convergence results of the present algorithm with respect to the spatial element numbers $M_x = M_y$ for various values of the direction parameter couples $(\theta_x, \theta_y)$, $K = 6, N = 10, S = 10, D = 1$ and $t_f = 0.5$.

**Problem 4 [11]**

Consider the (3+1) dimensional advection-diffusion equation (25) with $V_x = V_y = V_z = D_x = D_y = D_z = 1$ for which the exact solution is given by [11]

$$u(x,y,z,t) = \frac{1}{(4t+1)^{3/2}} exp\left(-\frac{(x-t-0.5)^2}{4t+1} - \frac{(y-t-0.5)^2}{4t+1} - \frac{(x-t-0.5)^2}{4t+1}\right), \ t > 0 \text{ and } (x,y,z) \in$$

$[0, 1] \times [0, 1] \times [0, 1].$  (51)

The nonhomogeneous time-dependent Dirichlet boundary conditions can be taken from exact solution (42), and the initial condition is considered as

$u(x,y,z,0) = exp(-(x-0.5)^2 - (y-0.5)^2 - (z-0.5)^2).$  (52)

The order refinement procedure is more vital in three-dimensional equations to minimize the computational cost of the algorithms. It is proven that the forward, backward, and central directional IELDM provides exponential convergence with increasing order of the method in Figure 9. The parameter values $M_x = M_y = M_z = 2, N = 10, S = 4$ and $t_f = 0.1$ are used to produce the figure with increasing orders from $K = 5$ to $K = 10$. The $h$-refinement procedure is known to be another way of increasing accuracy with the use of more elements in computational domains. The $h$-refinement response of the IELDTM is illustrated in Figure 10 for various continuity directions and the parameter values $K = 6, N = 10, S = 4$ and $t_f = 0.5$. As observed from the figure, the IELDTM ensures power convergence with increasing values of spatial elements. The IELDTM solution of Problem 4 and the corresponding pointwise absolute errors are illustrated with the consideration of the parameter values $K = 12, M_x = M_y = M_z = 1, N = 10, S = 4$ and $t_f = 0.25$ in Figure 10. As observed from the presented contour plots, the IELDTM accurately captures the physical behaviour using optimum degrees of freedom. As proven in Problem 2, the IELDTM has the ability to catch accurate physical behaviour of the advection-dominated problems in 2D. Since we perform the simulations on our personal computers, the advection-dominated 3D cases lead to a storage problem. It is believed that the IELDTM can solve these challenging problems with the optimized degrees of freedom with the computers that have relatively large RAM capacity.

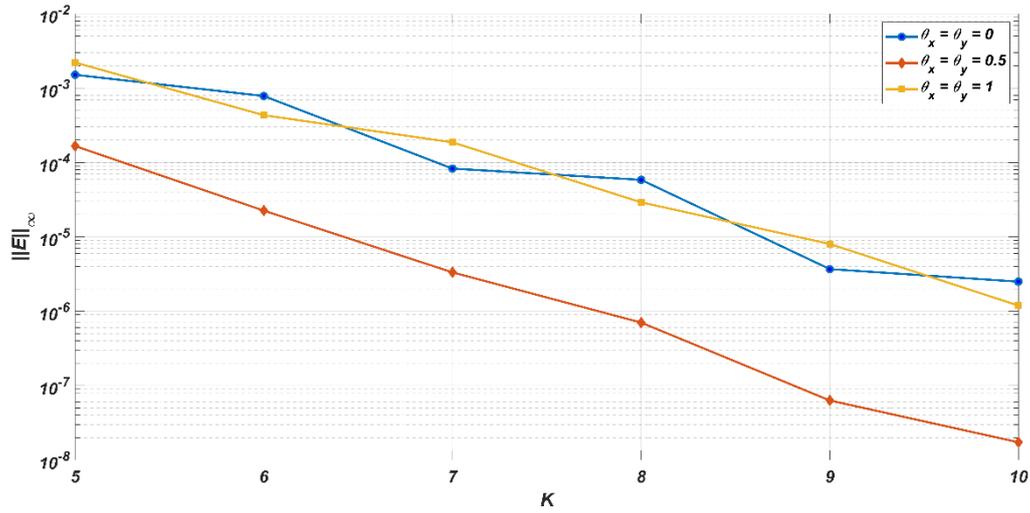

**Figure 9:** Convergence results of the present algorithm with respect to spectral convergence order $K$ for various values of the direction parameter triples $(\theta_x, \theta_y, \theta_z)$, $M_x = M_y = M_z = 2, N = 10, S = 4$ and $t_f = 0.1$.

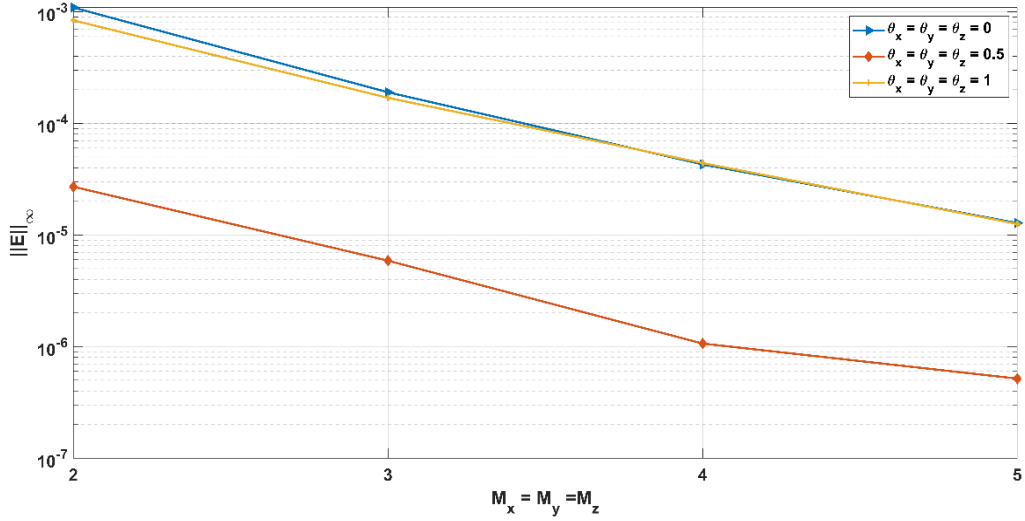

**Figure 10:** Convergence results of the present algorithm with respect to the spatial element numbers $M_x = M_y = M_z$ for various values of the direction parameter triples $(\theta_x, \theta_y, \theta_z)$, $K = 6, N = 10, S = 4$ and $t_f = 0.5$.

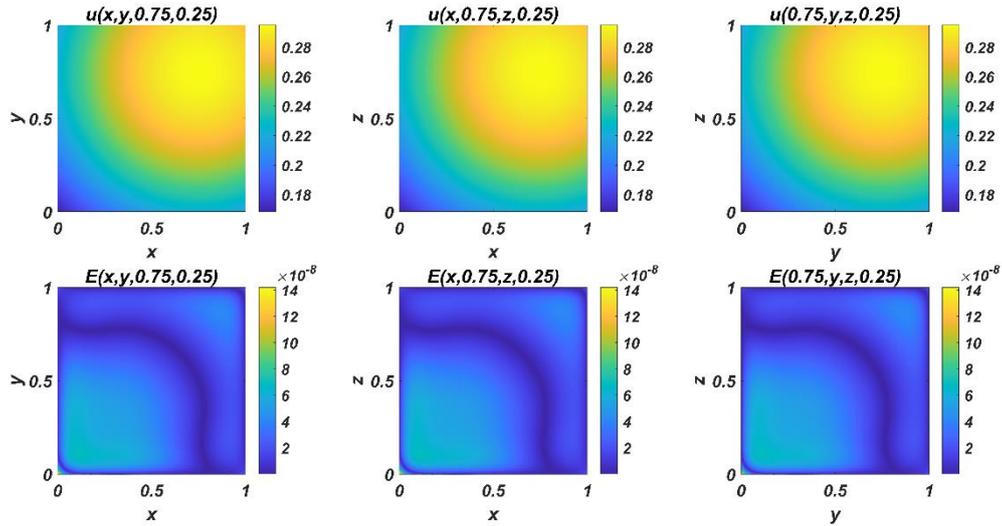

**Figure 11:** Numerical results and corresponding pointwise absolute errors of Problem 4 with the use of the IELDTM for $K = 12, M_x = M_y = M_z = 1, N = 10, S = 4$ and $t_f = 0.25$.

## 5 Conclusions and Recommendations

In this study, an implicit-explicit local differential transform method (IELDTM) has been derived for the (2+1) and (3+1) dimensional advection-diffusion equations. The parabolic PDEs have been reduced to a system of nonhomogeneous elliptic PDEs with the Chebyshev spectral collocation method. The present IELDTM has been proven to be a higher-order, stability-preserved, and versatile numerical technique for solving higher dimensional AD equations. The experimental converge analysis has shown that the IELDTM provides excellent convergence properties with both $h$−refinement and $p$−refinement. The IELDTM leads to the reduction of the required degrees of freedom for achieving higher-order convergent numerical results. It has been shown that the IELDTM needs only $(2K+1)$ unknowns to represent a complete two-dimensional polynomial of order $K$ while the entire polynomial includes $\frac{(K+1)(K+2)}{2}$ terms. Similarly, the IELDTM needs only $(K+1)^2$ unknowns to represent a complete three-dimensional polynomial of order $K$ while the entire polynomial includes $\frac{(K+1)(K+2)(K+3)}{6}$ terms. By comparing with the existing finite difference methods, finite element methods, and spectral methods, the IELDTM has been shown to supply more accurate results with far fewer degrees of freedom than both methods.

**Acknowledgements**
The first author would like to thank the Science Fellowships and Grant Programs Department of TUBITAK (BIDEB) for their support to his academic research.

**References**

[1] W. HUNDSDORFER AND J. VERWER, *Numerical solution of time-dependent advection-diffusion-reaction equations*, Springer-Verlag, 2003.

[2] R.C. MITTAL AND A. TRIPATHI, *Numerical solutions of two-dimensional Burgers' equations using modified Bi-cubic B-spline finite elements*, Eng. Comput. 32(5) (2015), pp. 1275-1306.

[3] D. FENGA, I. NEUWEILERA, U. NACKENHORST AND T. WICK, *A time-space flux-corrected transport finite element formulation for solving multi-dimensional advection-diffusion-reaction equations*, J. Comput. Phys., 396 (2019), pp. 31–53.

[4] X. LIU, J. WANG AND Y. ZHOU, *A space–time fully decoupled wavelet Galerkin method for solving two-dimensional Burgers' equations*, Comput. Math. with Appl., 72(12) (2016), pp. 2908-2919.

[5] E. TOMBAREVIĆ, V.R. VOLLER AND I. VUŠANOVIĆ, Detailed CVFEM Algorithm for Three Dimensional Advection-diffusion Problems, Comput. Model. Eng. Sci., 96(1) (2013), pp. 1–29.

[6] J. ESHAGHI, S. KAZEM AND H. ADIBI, *The local discontinuous Galerkin method for 2D nonlinear time-fractional advection–diffusion equations*, Eng. Comput., 35 (2019), pp. 1317–1332.

[7] L.T. DIOSADYA AND S.M. MURMAN, *Scalable tensor-product preconditioners for high-order finite-element methods: Scalar equations*, J. Comput. Phys., 394 (2019), pp. 759-776.

[8] N. MISHRA AND YVSS. SANYASIRAJU, *Exponential compact higher-order schemes and their stability analysis for unsteady convection-diffusion equations*, Int. J. Comput. Methods, 11(1) (2014), pp. 1350053.

[9] J.I. RAMOS, *A conservative method of lines for advection-reaction-diffusion equations*, Int. J. Numer. Method. H., 30(11) (2020), pp. 4735-4763.

[10] S. SINGH, D. BANSAL, G. KAUR AND S. SIRCAR, *Implicit-explicit-compact methods for advection diffusion reaction equations*, Comput. Fluids, 212 (2020), pp. 104709.

[11] S. KARAA, *A high-order compact ADI method for solving three-dimensional unsteady convection-diffusion problems*, Numer. Methods Partial Differ. Equ., 22(4) (2006), pp. 983-993.

[12] M. MANZINI AND A. RUSSO, A finite volume method for advection–diffusion problems in convection-dominated regimes, Comput. Methods Appl. Mech. Engrg., 197 (2008), pp. 1242–1261.

[13] J. SCHÄFER, X. HUANG, S. KOPECZ, P. BIRKEN, M.K. GOBBERT AND A. MEISTER, *A Memory-efficient finite volume method for advection-diffusion-reaction systems with non-smooth sources*, Numer. Methods Partial Differ. Equ., 31(1) (2015), pp. 143-167.

[14] M. HUSSAIN AND S. HAQ, *Weighted meshless spectral method for the solutions of multi-term time fractional advection-diffusion problems arising in heat and mass transfer*, Int. J. Heat Mass Tran., 129 (2019), pp. 1305–1316.

[15] A. KARAMANLI AND A. MUGAN, *Strong form Meshless Implementation of Taylor Series Method*, Appl. Math. Comput., 219 (2013), pp. 9069–9080.

[16] F. YAN, J.H. LV, X.T. FENG AND P.Z. PAN, *A new hybrid boundary node method based on Taylor expansion and the Shepard interpolation method*, Int. J. Numer. Meth. Engng. 102 (2015), pp. 1488–1506.

[17] J. YANGA, H. HUA, Y. KOUTSAWA AND M. POTIER-FERRY, *Taylor meshless method for solving nonlinear partial differential equations,* J. Comput. Phys., 348 (2017), pp. 385–400.

[18] J. YANG, H. HUA AND M. POTIER-FERRY, *Computing singular solutions to partial differential equations by Taylor series*, C. R. Mecanique, 346 (2018), pp. 603–614.

[19] J. YANG, Y. KOUTSAWA, M. POTIER-FERRY AND H. HU, Changing variables in Taylor series with applications to PDEs, Eng. Anal. Bound. Elem., 112 (2020), pp. 77–86.

[20] A. KARAMANLI, *Radial basis Taylor series method and its applications*, Eng. Comput. Vol. ahead-of-print No. ahead-of-print. https://doi.org/10.1108/EC-05-2020-0256

[21] R. ABAZARI AND A. BORHANIFAR, *Numerical study of the solution of the Burgers and coupled Burgers equations by a differential transformation method*, Comput. Math. with App., 59(8) (2010), pp. 2711-2722.

[22] A.H. KHATER, R.S. TEMSAH AND M.M. HASSAN, *A Chebyshev spectral collocation method for solving Burgers'-type equations*, J. Comput. Appl. Math., 222 (2008), pp. 333–350.

[23] Z.F. TIAN AND Y.B. GE, A fourth-order compact ADI method for solving two dimensional unsteady convection diffusion problems, J. Comput. Appl. Math., 198(1) (2007), pp. 168–286.